\RequirePackage[english]{babel}
\documentclass[a4paper,twoside]{article}

\usepackage{amssymb, amsfonts, amsmath, amsthm, epsfig, amscd, wrapfig}
\usepackage{mathrsfs, enumerate, caption, fancyhdr, tabularx}
\usepackage{booktabs, vmargin, multirow, rotating}
\usepackage[vcentermath]{youngtab}

\usepackage[all]{xy}
\usepackage[english]{babel}
\usepackage[latin1]{inputenc}

\theoremstyle{definition}
\newtheorem{lemma}{Lemma}

\newtheorem{teo}[lemma]{Theorem}
\newtheorem{prop}[lemma]{Proposition}
\newtheorem{cor}[lemma]{Corollary} 

\theoremstyle{remark}
\newtheorem{oss}{Remark}

\newcommand{\res}[2]{\ensuremath \textrm{Res}_{#1}^{#2}}
\newcommand{\ind}[2]{\ensuremath \textrm{Ind}_{#1}^{#2}}

\newcommand{\scr}[1]{\ensuremath {\mathscr{#1}}}

\newcommand{\sx} {\ensuremath {\left(}}
\newcommand{\dx} {\ensuremath {\right)}}
\newcommand{\mat}[1]{\ensuremath {\mathbb{#1}}}

\newcommand{\gen}[1]{\ensuremath \langle #1 \rangle}
\newcommand{\agisce} {\ensuremath \curvearrowright}

\renewcommand{\tilde}[1]{\ensuremath \widetilde{#1}}
\renewcommand{\epsilon}{\ensuremath \varepsilon}

\newcommand{\smalldiag}{\ensuremath \Yboxdim4.5pt}

\title{Symmetric group actions on the cohomology of configurations in $\mat R^d$}

\author{Giacomo d'Antonio
  and Giovanni Gaiffi
  }

\begin{document}
\maketitle

\begin{abstract}
  In this paper we deal with  the action of the symmetric group  on the cohomology of the configuration space $C_n(d)$ of
  $n$ points in $\mat R^d$. This topic has been studied by several authors, e.g. \cite{cohen1993representation},
  \cite{lehrer1986asg}, \cite{lehrer1987psa}, \cite{lehrer2000ecc}, \cite{mathieu1996hsigman} and \cite{gaiffi1996asn}.
  It is well-known that, for $d$ even, $H^*(C_n(d); \mat C) \cong 2\ind{S_2}{S_n}1$ and, for $d$ odd,
  $H^*(C_n(d); \mat C) \cong \mat C S_n$.

  On the cohomology algebra $H^*(C_n(d); \mat C)$ there is, in addition to the natural $S_n$-action, an \emph{extended action} 
  of $S_{n+1}$; this was first shown for the case when $d$ is even in \cite{mathieu1996hsigman}, \cite{gaiffi1996asn} and \cite{robinson1996trs}.
  For the case when $d$ is odd it was shown in \cite{mathieu1996hsigman} (anyway we will give an elementary algebraic construction
  of the extended action for this case).
  
  The purpose of this article is to present 
  several results that can be obtained,
  in an elementary way,
  exploiting the interplay between the extended
  action and the standard action. Among these we will recall a quick proof for the formula cited above for the case when $d$
  is even and show how to extend this proof to the case when $d$ is odd.
  We will also show how to locate among the homogeneous components of the graded algebra $H^*(C_n(d); \mat C)$ the copies
  of the standard, sign and standard tensor sign representations and we will give  explicit formulas for both the
  extended and the canonical actions on the low-degree cohomology modules.
\end{abstract}

\section{Introduction}
\label{sec:introduction}
We are concerned with the action of the symmetric group $S_n$ on the
cohomology algebra of the \emph{configuration space of $n$ points
in $\mat R^d$}, that is the space
\begin{displaymath}
  C_n(d) = \left\{(p_1, \dots, p_n) \in \sx\mat R^d\dx^n:\, p_i \neq
  p_j,\,\forall i \neq j\right\}.
\end{displaymath}
$S_n$ acts on $C_n(d)$ permuting coordinates and this action
induces an action on the cohomology algebra $H^*(C_n(d); \mat C)$.

In the particular case $d = 2$, $C_n(2) = M(\scr B_{n})$ is the
complement of the (complex)  braid arrangement  \(\scr B_{n}\). 
This a widely studied object,
in \cite{arnold1969crc} Arnol'd gave a presentation for the cohomology
algebra $H^*(C_n(2); \mat C)$; he proved that it is the skew-commutative
algebra with generators $\{A_{i,j}:\, i < j\}$ of degree $1$ and
relations:
\begin{equation}
  \label{eq:arnold}
  A_{i,j}A_{i,k} - A_{i,j}A_{j,k} + A_{i,k}A_{j,k} = 0,\quad
  \forall i < j < k.
\end{equation}
This result can be generalized to hyperplane 
arrangements (cfr. \cite{orlik1980cat}, \cite{orlik1992ah},
\cite{yuzvinsky2001osa} and \cite{bjorner1992csc}).

The action of $\sigma \in S_n$ on $H^*(C_n(2); \mat C)$ can be described as $\sigma A_{i,j} =
A_{\sigma i, \sigma j}$ and was  studied
by Lehrer and Solomon in \cite{lehrer1986asg} and by Lehrer in \cite{lehrer1987psa}.
Among other results they proved a formula for the character
of this action; precisely the following holds:
\begin{equation}\label{eq:character2}
  \chi_{H^*(C_n(2);\mat C)} = 2 \ind{S_2}{S_n}(1)
\end{equation}
where $1$ is the character of the trivial representation of $S_2$.
In \cite{gaiffi1996asn} the second author  gave a quick   proof for this formula
introducing an action of $S_{n+1}$ on $H^*(C_n(2);\mat C)$ which
restricts to the natural $S_n$ action (this action
is called the \emph{extended} or \emph{hidden}  action). 
A similar approach was exploited by
Mathieu in \cite{mathieu1996hsigman} and by Robinson and Whitehouse in
\cite{robinson1996trs}.

In the general case the algebra $H^*(C_n(d); \mat C)$ can be presented as
follows (cfr. \cite{cohen1976hil} and \cite{cohen1993representation}): it is the associative graded
algebra with generators $\{A_{i,j}: 1 \leq i,j \leq n\}$ 
(with \(i\neq j\); for convenience of notation we allow  \(i=j\) and then  \(A_{i,i}=0\)) of
degree $d-1$ and relations
\begin{eqnarray}
& & \label{eq:rel1} A_{i,j} = (-1)^d A_{j,i},\\
& & \label{eq:rel2} A_{i,j}A_{h,k} = (-1)^{d-1}A_{h,k}A_{i,j},\\
& & \label{eq:rel3} A_{i,j}A_{i,k} = A_{k,j}(A_{i,k} - A_{i,j}) 
\mbox{ for } 1 \leq i \leq n \mbox{ and } j \leq k.
\end{eqnarray}
Again these results generalize to the complement
of a subspace arrangement (cfr. \cite{feichtner2000cac},
\cite{delongueville2001crc} and \cite{deligne2000algebre}).

The action of $S_n$ on $H^*(C_n(d); \mat C)$ for arbitrary $d$ 
was studied by Cohen and Taylor in \cite{cohen1993representation}
and by Lehrer in \cite{lehrer2000ecc}.    
Lehrer provided formulas for the generalized Poincar\'e polynomials associated to the representations $H^*(C_n(d); \mat C)$. 
It turns out that there is a qualitative
difference between the case when $d$ is even and the case when
$d$ is odd. The argument for the case $d = 2$ translates literally to
the case when $d$ is even (cfr. \cite{lehrer2000ecc}) and formula
\eqref{eq:character2} still holds. Also the construction of the
extended action can be translated to the case when $d$ is even.
For the case when $d$ is odd both Cohen and Taylor in \cite{cohen1993representation}
and Lehrer in \cite{lehrer2000ecc} proved, with different arguments, 
that $H^*(C_n(d); \mat C)$ is the regular representation $\mat C S_n$.

We will construct (section \ref{sec:extended}) 
in an elementary  way  an $S_{n+1}$ action on
$H^*(C_n(d);\mat C)$ for the case when $d$ is odd (the same action was described with a different method
in \cite{mathieu1996hsigman}) and use it in section \ref{sec:character}  to prove 
quickly some results of \cite{cohen1993representation} and \cite{lehrer2000ecc}.
In addition we will show how the extended action can be used, both in the case when $d$ is
even and in the case when $d$ is odd, to locate 
the copies of the standard, sign and standard tensor sign representations
of $S_n$ on the homogeneous components $H^{k(d-1)}(C_n(d); \mat C)$ (section \ref{sec:location}) and
to prove explicit formulas for the decomposition of the degrees
$d-1$ and $2(d-1)$ (section \ref{sec:lowdegrees}).

\section{The extended $S_{n+1}$ action}
\label{sec:extended}
We now discuss the definition of an  extended action on $H^*(C_n(d); \mat C)$. We distinguish
the case when $d$ is odd and the case when $d$ is even.
In the former case we see from relations \eqref{eq:rel1}-\eqref{eq:rel3}
that there is an isomorphism of graded $S_n$-modules
\begin{equation}
  \label{eq:isodeven}
  H^*(C_n(d); \mat C) \to H^*(C_n(2); \mat C) \otimes 1
\end{equation}
where $1$ is the graded $S_n$-module whose only non-zero component
is the trivial representation at degree $d-1$. There are (at least)
three different ways to extend the action of $S_n$ on
$H^*(C_n(2); \mat C)$ to an $S_{n+1}$-action (see \cite{gaiffi1996asn},
\cite{mathieu1996hsigman} and \cite{robinson1996trs});
the isomorphism \eqref{eq:isodeven} allows us to 
carry this extended action to $H^*(C_n(d); \mat C)$.

In the case when $d$ is odd we can rewrite the relations
\eqref{eq:rel1}-\eqref{eq:rel3} as follows:
\begin{eqnarray*}
& & A_{i,j} = -A_{j,i},\\
& & A_{i,j}A_{h,k} = A_{h,k}A_{i,j},\\
& & A_{i,j}A_{i,k} - A_{i,j}A_{j,k} + A_{i,k}A_{j,k} = 0
\end{eqnarray*}

We first look at the degree $d-1$; let $V = \mat C^n$ be the permutation representation,
we have an equivariant isomorphism of \(S_n\) modules
\begin{displaymath}
  \begin{array}{c}
    \bigwedge^2 V \to H^{d-1}(C_n(d); \mat C)\\
    e_i \wedge e_j \mapsto A_{i,j}.
  \end{array}
\end{displaymath}

The action of $S_n$ on $V$ can be extended to an $S_{n+1}$-action;
from Pieri's rule we see that any extended action must
be isomorphic to the standard representation of $S_{n+1}$,
that is $V_{(n,1)} = \ker (x_0 + \cdots + x_n) \subseteq \mat C^{n+1}$.
We choose a basis for $V_{(n,1)}$ of elements
$\{v_1,\dots, v_n\}$ where $v_i = e_i - e_0$; identifying
$S_n = \{\sigma \in S_{n+1}:\, \sigma(0) = 0\}$ we have an $S_n$-equivariant
isomorphism
\begin{equation}\label{eq:isostandard}
  \begin{array}{c}
    \res{S_n}{S_{n+1}} V_{(n,1)} \to V\\
    v_j \mapsto e_j
  \end{array}
\end{equation}
and we can define the $S_{n+1}$ action on $V$ as the unique action
that makes \eqref{eq:isostandard} into an equivariant  isomorphism
$V_{(n,1)} \to V$.

The $S_{n+1}$ action on $V$ induces an $S_{n+1}$ action on $\bigwedge^2 V \cong
H^{d-1}(C_n(d);\mat C)$. We can describe this action as follows: 
if $\sigma \in S_n$ then $\sigma A_{i,j} = A_{\sigma i, \sigma j}$ and
\begin{eqnarray}
  \label{eq:s0action1}
  (0,1)A_{i,j} = A_{i,j} - A_{1,j} + A_{1,i} & \mbox{ if } 1 < i < j\\
  \label{eq:s0action2}
  (0,1)A_{1,j} = -A_{1,j}.
\end{eqnarray}

In particular we are able to decompose $H^{d-1}(C_n(d); \mat C)$
for every $n$ and $d$ odd:
\begin{prop}
  \label{prop:decompositionhd-1odd}
  There is an isomorphism of $S_{n+1}$-modules
  \begin{displaymath}
    H^{d-1}(C_n(d); \mat C) \cong \bigwedge^2 V_{(n,1)} \cong V_{(n-1,1,1)}.
  \end{displaymath}
  Using Pieri's rules we obtain the following isomorphism of
  $S_n$-modules
  \begin{equation}
    \label{eq:decompositionhd-1odd}
    H^{d-1}(C_n(d); \mat C) = V_{(n-1,1)} \oplus V_{(n-2,1,1)}.
  \end{equation}
\end{prop}

From relations \eqref{eq:rel1}-\eqref{eq:rel3} we see that there is
an equivariant isomorphism
\begin{displaymath}
  H^*(C_n(d);\mat C) \cong S(H^{d-1}(C_n(d);\mat C))/I_{n,d}
\end{displaymath}
where $S(H^{d-1}(C_n(d);\mat C))$ is the symmetric algebra
on $H^{d-1}(C_n(d);\mat C)$ and $I_{n,d}$ is the ideal of
relations:
\begin{displaymath}
  I_{n,d} = \gen{A_{i,j}A_{i,k} - A_{k,j}(A_{i,k} - A_{i,j}):
    \mbox{ for } 1 \leq i \leq n \mbox{ and } j \leq k}.
\end{displaymath}
In particular, in order to extend the $S_n$ action
on $H^*(C_n(d); \mat C)$ we only need to prove
that the ideal $I_{n,d}$ is invariant under the action
of $S_{n+1}$ on $S(H^{d-1}(C_n(d);\mat C))$.
This is indeed the case and is proved with a short  explicit
computation. One has to check the equalities (in $H^*(C_n(d); \mat C)$)
\begin{displaymath}
  (0,1) A_{i,j}A_{i,k} = (0,1) A_{k,j}(A_{i,k} - A_{i,j}).
\end{displaymath}
We notice that, since the expression above is symmetric in \(j\) and \(k\),  it suffices to  distinguish three cases: the case when $i = 1$,
the case when $j = 1$ and the case when $i,j,k \neq 1$.

\section{The character of the $S_n$ action on $H^*(C_n(d); \mat C)$}
\label{sec:character}
In this section we will recall some results and proofs from \cite{gaiffi1996asn}, 
\cite{mathieu1996hsigman} and \cite{robinson1996trs} and 
use them  to show a quick proof of formula
$H^*(C_n(d); \mat C) \cong \mat C S_n$ for $d$ odd.

When $d$ is even isomorphism \eqref{eq:isodeven} provides us an analogous of
\cite[Theorem 4.1]{gaiffi1996asn} (see also
\cite{mathieu1996hsigman} and \cite{robinson1996trs}), i.e.
\begin{equation}
  \label{eq:recursivedeven}
  H^{k(d-1)}(C_n(d); \mat C) \cong H^{k(d-1)}(C_{n-1}(d); \mat C) \oplus
  \sx H^{(k-1)(d-1)}(C_{n-1}(d); \mat C) \otimes V_{(n-1,1)}\dx
\end{equation}
which connects the canonical $S_n$-action (on the left) with
the extended $S_n$-action on $H^*(C_{n-1}(d); \mat C)$ (on the right).

Consider now the case when $d$ is odd; let $\eta: H^*(C_{n-1}(d);\mat C) \to H^*(C_n(d); \mat C)$
be the map $A_{i,j} \mapsto A_{i,j}$ (i.e. the map induced by the projection on the
first $n-1$ factors $C_n(d) \to C_{n-1}(d)$). If we call
$s_j = (j,j+1) \in S_{n+1}$ we have from formulas
\eqref{eq:s0action1} and \eqref{eq:s0action2}
that the map $\eta$ is $\gen{s_0, \dots, s_{n-2}}$-equivariant.
Recall the following well known result (see \cite{cohen1993representation}):
\begin{prop}\label{prop:admissiblebasis}
  The algebra $H^*(C_n(d); \mat C)$ has a basis given by the  elements
  \begin{displaymath}
    A_{i_1,j_1}A_{i_2,j_2} \cdots A_{i_k,j_k}
  \end{displaymath}
  with $i_h < j_h$ and $1 < j_1 < j_2 < \cdots < j_k \leq n$.
\end{prop}
Such elements are usually called \emph{admissible monomials}.
We will write $\chi(n,k)$ for the character of the action of $S_n$ on
$H^k(C_n(d); \mat C)$ and $\tilde\chi(n,k)$ for the character of the
extended action of $S_{n+1}$ on $H^k(C_n(d); \mat C)$.

With these ingredients we can translate almost verbatim the proof
of \cite[Theorem 4.1]{gaiffi1996asn} to obtain the following
result, which we state for arbitrary $d$ 
(see also \cite[Theorem 4.4]{mathieu1996hsigman}).
\begin{teo}
  \label{teo:recursiveformula}
  For any $n,k,d$ it holds:
  \begin{displaymath}
    \chi(n,k) = \tilde\chi(n-1,k) + p_n\tilde\chi(n-1,k-1),    
  \end{displaymath}
  where $p_n$ is the character of the standard representation of $S_n$.  
\end{teo}
\begin{proof}
  We discuss only the case when $d$ is odd.
  Consider the $\gen{s_0, \dots, s_{n-2}}$-submodule
  $\Omega_{n-1,k} = \eta\sx H^k(C_{n-1}(d); \mat C)\dx \subseteq
  H^k(C_n(d);\mat C)$. We can write
  \begin{equation}
    \label{eq:decomposition1}
    H^k(C_n(d); \mat C) = \Omega_{n-1,k} \oplus N\cdot\Omega_{n-1,k-1}
  \end{equation}
  where $N = \oplus_{j = 1}^{n-1} \mat C A_{j,n}$ is certainly
  $\gen{s_1, \dots, s_{n-2}}$-invariant but, in general,
  is not an $\gen{s_0, \dots, s_{n-2}}$-submodule.

  Now consider the case $k = 1$, we have
  \begin{displaymath}
    H^1(C_n(d); \mat C) = \Omega_{n-1,1} \oplus N = \Omega_{n-1,1} \oplus T
  \end{displaymath}
  where $T$ is an $\gen{s_0, \dots, s_{n-2}}$-invariant complement of
  $\Omega_{n-1,1}$ (in particular its restriction to
  $\gen{s_1, \dots, s_{n-2}}$ is isomorphic to $N$).
  But $S_{n-1} = \gen{s_1, \dots, s_{n-2}}$ permutes the elements
  $A_{1,n}, \cdots, A_{n-1,n}$ and therefore
  \begin{displaymath}
    \res{S_{n-1}}{S_n} T \cong N \cong V_{(n-1)} \oplus V_{(n-2,1)}
  \end{displaymath}
  where $V_{(n-1)}$ is the trivial representation and
  $V_{(n-2,1)}$ is the standard representation. By Pieri's rule
  we have $T \cong V_{(n-1,1)}$ as $S_n$-module.

  We can still write
  \begin{displaymath}
    H^k(C_n(d); \mat C) = \Omega_{n-1,k} \oplus T\cdot\Omega_{n-1,k-1},
  \end{displaymath}
  indeed we have $A_{i,n} \in H^1(C_n(d);\mat C) \Rightarrow
  A_{i,n} = \gamma_i^{(1)} + \gamma_i^{(2)}$ with
  $\gamma_i^{(1)} \in \Omega_{n-1,1}$ and $\gamma_i^{(2)} \in T$.
  Now let $z = z_0 + \sum_{j = 1}^n A_{j,n} z_j \in H^k(C_n(d); \mat C)$ 
  with $z_0 \in \Omega_{n-1,k}$ and for $j > 0$, $z_j \in \Omega_{n-1,k-1}$; 
  then we have
  \begin{displaymath}
    z = \underbrace{z_0 + \sum_{j = 1}^n \gamma_j^{(1)}z_j}_{\in \Omega_{n-1,k}} + 
    \underbrace{\sum_{j=1}^n \gamma_j^{(2)} z_j}_{\in T\cdot\Omega_{n-1,k-1}}.
  \end{displaymath}
  Therefore $H^k(C_n(d); \mat C) = \Omega_{n-1,k} + T\cdot\Omega_{n-1,k-1}$ and
  the sum is direct by a dimension argument.
  In particular we have 
  $\dim T\cdot\Omega_{n-1,k-1} = \dim (T \otimes \Omega_{n-1,k-1})$
  and therefore there is an equivariant isomorphism $T\cdot\Omega_{n-1,k-1} \cong
  T \otimes \Omega_{n-1,k-1}$. 

  We have proved a decomposition of
  $\gen{s_0, \dots, s_{n-2}}$-modules
  \begin{displaymath}
    H^k(C_n(d); \mat C) \cong
    H^k(C_{n-1}(d); \mat C) \oplus \sx P_n \otimes H^{k-1}(C_{n-1}(d); \mat C)\dx.
  \end{displaymath}
  Now consider the subgroups $H_1 = \gen{s_0, \dots, s_{n-2}}$ and
  $H_2 = \gen{s_1, \dots, s_{n-1}}$ of $S_{n+1}$; these are conjugate
  subgroups and therefore
  \begin{displaymath}
    \res{H_1}{S_{n+1}} H^k(C_n(d);\mat C) \cong \res{H_2}{S_{n+1}} H^k(C_n(d); \mat C)
  \end{displaymath}
  and the term on the right is the natural $S_n$ action on $H^k(C_n(d); \mat C)$.
\end{proof}

As a consequence we immediately have the following.
\begin{cor}
  \label{cor:inducedrepresentation}
  For any $n > 2$ and any \(d \geq 2\) the following equality of $S_n$-modules holds:
  \begin{displaymath}
    H^*(C_n(d); \mat C) = \ind{S_{n-1}}{S_n} H^*(C_{n-1}(d); \mat C).
  \end{displaymath}
\end{cor}
\begin{proof}
  Call $\chi_n = \sum_{k=0}^{n-1} \chi(n,k)$ the character of the
  action of $S_n$ on $H^*(C_n(d); \mat C)$ and
  $\tilde\chi_n = \sum_{k=0}^{n-1} \tilde\chi(n,k)$ the character of
  the extended action. Then from theorem \ref{teo:recursiveformula}
  and from the fact that $\tilde\chi(-1,n-1) = \tilde\chi(n-1,n-1) = 0$
  we have
  \begin{displaymath}
    \chi_n = \sum_{k = 0}^{n-1} \chi(n, k) = 
    \sum_{k=0}^{n-1}\sx\tilde\chi(n-1,k) + p_n\tilde\chi(n-1,k-1)\dx = (1 + p_n)\tilde\chi_{n-1}.
  \end{displaymath}
  Recall that if $H \subseteq G$ is a subgroup and $M$ is a $G$-module we have
  $\ind{H}{G}\res{H}{G} M = M \otimes \ind{H}{G}(1)$. In our case we have
  \begin{displaymath}
    \ind{S_{n-1}}{S_n} \chi_{n-1} = \ind{S_{n-1}}{S_n} \res{S_{n-1}}{S_n} \tilde\chi_{n-1} =
    (\ind{S_{n-1}}{S_n} 1)\tilde\chi_{n-1} = (1 + p_n)\tilde\chi_{n-1} = \chi_n
  \end{displaymath}  
\end{proof}

As remarked in \cite[Theorem 4.4]{gaiffi1996asn}, corollary \ref{cor:inducedrepresentation} provides a quick proof of Lehrer and Solomon
result for $d$ even: $H^*(C_n(d); \mat C) = 2\, \ind{S_2}{S_n}\,1$, since $H^*(C_2(d); \mat C)$ consists of two copies of the trivial representation of \(S_2\).
Analogously, when \(d\) is odd  we can now  prove the following result of \cite{cohen1993representation} and \cite{lehrer2000ecc}:
\begin{teo}
  \label{teo:regular}
  When $d$ is odd we have:
  \begin{displaymath}
    H^*(C_n(d); \mat C) \cong \mat C S_n.
  \end{displaymath}
\end{teo}
\begin{proof}
  By induction on $n$; it is easy to check that $H^*(C_2(d); \mat C) \cong \mat C S_2$ 
  (we have $H^0(C_2(d); \mat C) \cong {\tiny\yng(2)}$ and
  $H^{d-1}(C_2(d); \mat C) \cong {\tiny\yng(1,1)}$).
  Now, using the inductive hypothesis and corollary \ref{cor:inducedrepresentation}
  we have
  \begin{displaymath}
    H^*(C_n(d); \mat C) = \ind{S_{n-1}}{S_n} H^*(C_{n-1}(d); \mat C) \cong 
    \ind{S_{n-1}}{S_n} \mat C S_{n-1} \cong \mat C S_n.
  \end{displaymath}
\end{proof}

For low $n$, the recursive relation of Theorem \ref{teo:recursiveformula} allows us to compute the graded character of the \(S_n\) action, 
as is shown in tables \ref{table:decomposition} and \ref{table:decompositionodd}.

\begin{sidewaystable}[tbp]
  \centering
  \begin{tabular*}{0.75\textwidth}{@{\extracolsep{\fill}}*{6}{c}}
    \multicolumn{2}{c}{degrees} & $0$ & $1$ & $2$ & $3$ \\
    \toprule
    \multirow{2}{*}{$n = 2$} & can. & $\smalldiag{\yng(2)}$ \\
    \addlinespace
    & ext. & $\smalldiag{\yng(3)}$\\
    \addlinespace
    \midrule[0.3pt]
    \multirow{2}{*}{$n = 3$} & can. & $\smalldiag{\yng(3)}$ & 
    $\smalldiag{\yng(2,1)}$\\
    \addlinespace
    & ext. & $\smalldiag{\yng(4)}$ & $\smalldiag{\yng(2,2)}$\\
    \addlinespace
    \midrule[0.3pt]
    \multirow{2}{*}{$n = 4$} & can. & $\smalldiag{\yng(4)}$ &
    $\smalldiag{\yng(2,2)} \oplus {\yng(3,1)}$ &
    $\smalldiag{\yng(3,1)} \oplus {\yng(2,1,1)}$\\
    \addlinespace
    & ext. & $\smalldiag{\yng(5)}$ & 
    $\smalldiag{\yng(3,2)}$ & $\smalldiag{\yng(3,1,1)}$\\
    \midrule[0.3pt]
    \addlinespace
    $n = 5$ & can. &
    $\smalldiag{\yng(5)}$ & 
    $\smalldiag{\yng(3,2)} \oplus {\yng(4,1)}$ &
    $\smalldiag 2\,{\yng(3,1,1)} \oplus {\yng(3,2)} \oplus
    {\yng(4,1)} \oplus {\yng(2,2,1)}$ &
    $\smalldiag{\yng(4,1)} \oplus {\yng(2,1,1,1)} \oplus
    {\yng(3,1,1)} \oplus {\yng(3,2)} \oplus {\yng(2,2,1)}$
  \end{tabular*}
  \caption{Decomposition of $H^*(M(d\scr A_{n-1}); \mat C)$.}
  \label{table:decompositiondeconing}

  \vspace{2\baselineskip}
  \centering

  \begin{tabular*}{\textwidth}{@{\extracolsep{\fill}}*{7}{c}}
    \multicolumn{2}{c}{degrees} & $0$ & $d-1$ & $2(d-1)$ & $3(d-1)$ & $4(d-1)$ \\
    \toprule
    \multirow{2}{*}{$n = 2$} & can. & $\smalldiag{\yng(2)}$
    & $\smalldiag{\yng(2)}$ \\
    \addlinespace
    & ext. & $\smalldiag{\yng(3)}$ & $\smalldiag{\yng(3)}$\\
    \addlinespace
    \midrule[0.3pt]
    \multirow{2}{*}{$n = 3$} & can. & $\smalldiag{\yng(3)}$ & 
    $\smalldiag{\yng(3)} \oplus {\yng(2,1)}$ & $\smalldiag{\yng(2,1)}$ \\
    \addlinespace
    & ext. & $\smalldiag{\yng(4)}$ & 
    $\smalldiag{\yng(4)} \oplus {\yng(2,2)}$ & $\smalldiag{\yng(2,2)}$ \\
    \addlinespace
    \midrule[0.3pt]
    \multirow{2}{*}{$n = 4$} & can. & $\smalldiag{\yng(4)}$ &
    $\smalldiag{\yng(4)} \oplus {\yng(2,2)} \oplus {\yng(3,1)}$ &
    $\smalldiag{\yng(2,2)} \oplus 2\,{\yng(3,1)} \oplus {\yng(2,1,1)}$ &
    $\smalldiag{\yng(3,1)} \oplus {\yng(2,1,1)}$\\
    \addlinespace
    & ext. & $\smalldiag{\yng(5)}$ & 
    $\smalldiag{\yng(5)} \oplus {\yng(3,2)}$ &
    $\smalldiag{\yng(3,2)} \oplus {\yng(3,1,1)}$ &
    $\smalldiag{\yng(3,1,1)}$\\
    \midrule[0.3pt]
    \addlinespace
    $n = 5$ & can. &
    $\smalldiag{\yng(5)}$ & 
    $\smalldiag{\yng(5)} \oplus {\yng(3,2)} \oplus {\yng(4,1)}$ &
    $\smalldiag 2\,{\yng(3,1,1)} \oplus 2\,{\yng(3,2)} \oplus
    2\,{\yng(4,1)} \oplus {\yng(2,2,1)}$ &
    $\smalldiag 2\,{\yng(4,1)} \oplus {\yng(2,1,1,1)} \oplus
    3\,{\yng(3,1,1)} \oplus 2\,{\yng(3,2)} \oplus 2\,{\yng(2,2,1)}$ &  
    $\smalldiag{\yng(4,1)} \oplus {\yng(2,1,1,1)} \oplus
    {\yng(3,1,1)} \oplus {\yng(3,2)} \oplus
    {\yng(2,2,1)}$
  \end{tabular*}
\caption{Decomposition of $H^*(C_n(d); \mat C)$ when $d$ is even.}
\label{table:decomposition}
\end{sidewaystable}

\begin{sidewaystable}[tbp]
  \centering
  \begin{tabular*}{\textwidth}{@{\extracolsep{\fill}}*{7}{c}}
    \multicolumn{2}{c}{degrees} & $0$ & $d-1$ & $2(d-1)$ & $3(d-1)$ & $4(d-1)$\\
    \toprule
    \multirow{2}{*}{$n = 2$} & can. & $\smalldiag{\yng(2)}$
    & $\smalldiag{\yng(1,1)}$ \\
    \addlinespace
    & ext. & $\smalldiag{\yng(3)}$ & $\smalldiag{\yng(1,1,1)}$\\
    \addlinespace
    \midrule[0.3pt]
    \multirow{2}{*}{$n = 3$} & can. & $\smalldiag{\yng(3)}$ & 
    $\smalldiag{\yng(2,1)} \oplus {\yng(1,1,1)}$ &
    $\smalldiag{\yng(2,1)}$\\
    \addlinespace
    & ext. & $\smalldiag{\yng(4)}$ & 
    $\smalldiag{\yng(2,1,1)}$ & $\smalldiag{\yng(2,2)}$\\
    \addlinespace
    \midrule[0.3pt]
    \multirow{2}{*}{$n = 4$} & can. & $\smalldiag{\yng(4)}$ &
    $\smalldiag{\yng(2,1,1)} \oplus {\yng(3,1)}$ &
    $\smalldiag{\yng(1,1,1,1)} \oplus {\yng(3,1)} \oplus {\yng(2,1,1)}
    \oplus 2\,{\yng(2,2)}$ &
    $\smalldiag{\yng(3,1)} \oplus {\yng(2,1,1)}$\\
    \addlinespace
    & ext. & $\smalldiag{\yng(5)}$ & $\smalldiag{\yng(3,1,1)}$ &
    $\smalldiag{\yng(1,1,1,1,1)} \oplus {\yng(3,2)} \oplus {\yng(2,2,1)}$ &
    $\smalldiag{\yng(3,1,1)}$\\
    \midrule[0.3pt]
    \addlinespace
    $n = 5$ & can. &
    $\smalldiag{\yng(5)}$ & 
    $\smalldiag{\yng(3,1,1)} \oplus {\yng(4,1)}$ &
    $\smalldiag{\yng(1,1,1,1,1)} \oplus {\yng(2,1,1,1)} \oplus
    {\yng(4,1)} \oplus 2\,{\yng(3,2)} \oplus 2\,\yng(2,2,1) \oplus
    {\yng(3,1,1)}$ &
    $\smalldiag 3\,{\yng(3,1,1)} \oplus 2\,{\yng(2,1,1,1)} \oplus
    {\yng(4,1)} \oplus 2\,{\yng(3,2)} \oplus 2\,{\yng(2,2,1)}$ & 
    $\smalldiag{\yng(4,1)} \oplus {\yng(2,1,1,1)} \oplus
    {\yng(3,1,1)} \oplus {\yng(3,2)} \oplus
    {\yng(2,2,1)}$
  \end{tabular*}
\caption{Decomposition of $H^*(C_n(d); \mat C)$ when $d$ is odd.}
\label{table:decompositionodd}
\end{sidewaystable}

\section{Locating some irreducible representations}
\label{sec:location}
Using the recursive formula of theorem \ref{teo:recursiveformula} it is possible to locate some irreducible representations of $S_n$
in the homogeneous components $H^{k(d-1)}(C_n(d); \mat C)$;
namely we will locate the copies of the 
\emph{standard}, the \emph{sign} and the \emph{standard tensor sign}
representations. As before we need to distinguish the case when
$d$ is even and the case when $d$ is odd.

\subsection{The case $d$ even}
\label{sec:locatedeven}
Using isomorphism \eqref{eq:isodeven} we reduce ourselves to
study the action of $S_n$ on $H^*(M(\scr B_{n}); \mat C)$; more
precisely we study the action of $S_n$ on the cohomology of the 
complement of the \emph{essential braid arrangement} $\scr A_{n-1}$
(i.e. the arrangement in $\mat C^n/\gen{(1,1,\dots,1)}$ induced
by $\scr B_n$ or equivalently the \emph{Coxeter arrangement} of type $A_{n-1}$). 

Recall the \emph{deconing construction} from the theory of arrangements;
i.e. the deconing of the 
essential braid arrangement is the arrangement
$d\scr A_{n-1}$ on the vector space $\mat C^{n-2}$ such that
$M(d\scr A_{n-1}) \cong M(\scr A_{n-1})/\mat C^*$.
There is an $S_{n+1}$-equivariant isomorphism of graded algebras 
(\cite[Proposition 2.2]{gaiffi1996asn})
\begin{equation}
  \label{eq:isodeconing}
  H^*(M(\scr A_{n-1}); \mat C) \cong H^*(M(d\scr A_{n-1}); \mat C) \otimes \mat C[\epsilon]/\epsilon^2
\end{equation}
where \(\epsilon \) has degree 1 and $S_{n+1}$ acts trivially on $\mat C[\epsilon]/\epsilon^2$.
Futhermore theorem \ref{teo:recursiveformula} and corollary \ref{cor:inducedrepresentation}
still hold for the $S_n$-module $H^*(M(d\scr A_{n-1}); \mat C)$.
There is an analogous of \eqref{eq:character2} for
$H^*(M(d\scr A_{n-1}); \mat C)$, namely:
\begin{equation}
  \label{eq:character3}
  H^*(M(d\scr A_{n-1}); \mat C) = \ind{S_2}{S_n} 1.
\end{equation}
Moreover isomorphism \eqref{eq:isodeconing} allows us to know
the location of an irreducible representation in $H^*(M(\scr A_{n-1}); \mat C)$
once we know its location in $H^*(M(d\scr A_{n-1}); \mat C)$. 
We recall  that a formula for the generalized Poincar\`e series associated to the \(S_{n+1}\) action on $H^*(M(d\scr A_{n-1}); \mat C)$ 
has been shown  in \cite{getzler1995ope}, given that \(M(d\scr A_{n-1})\) is homeomorphic to the moduli space \({\mathcal M}_{0,n+1}\) 
of genus zero \(n+1\)-pointed curves (and its minimal De Concini-Procesi wonderful model - see \cite{de1995wonderful} -  is isomorphic to the Deligne-Mumford compatification of \({\mathcal M}_{0,n+1}\)).

As before theorem \ref{teo:recursiveformula}  
suffices to compute the graded character of the $S_n$ action on $H^*(M(d\scr A_{n-1}), \mat C)$ for low $n$,
as is shown in table \ref{table:decompositiondeconing}.

As a first observation we see that formula \eqref{eq:character3}
and Fr\"obenius reciprocity allow us to know the number of copies
of each irreducible representation in the whole 
$H^*(M(d\scr A_{n-1}); \mat C)$; in particular
\begin{enumerate}[(i)]
\item there is only one copy of the trivial representation in $H^*(M(d\scr A_{n-1}); \mat C)$
  (and must be at the degree $0$),
\item there are $n-2$ copies of the standard representation in $H^*(M(d\scr A_{n-1}); \mat C)$,
\item there are no copies of the sign representation in $H^*(M(d\scr A_{n-1}); \mat C)$,
\item there is one copy of the standard tensor sign representation in $H^*(M(d\scr A_{n-1}); \mat C)$.
\end{enumerate}

We will use the notation $\chi^*(n,k)$ for the character of the action of $S_n$
on $H^k(M(d\scr A_{n-1}); \mat C)$ and $\tilde\chi^*(n,k)$ for the character
of the extended action of $S_{n+1}$ on $H^k(M(d\scr A_{n-1}); \mat C)$.
\begin{prop}\label{prop:wherearethestandard}
  For $n \geq 3$ there is exactly one copy of the standard representation $V_{(n-1,1)}$
  in $H^k(M(d\scr A_{n-1}); \mat C)$ for each $0 < k < n-1$.
\end{prop}
\begin{proof}
  By induction on $n$, for $n = 3,4,5$ it follows from an explicit computation
  (see table \ref{table:decompositiondeconing}). Let $n > 5$;  we have
  \begin{displaymath}
    \gen{\chi^*(n,k), p_n} = \gen{\tilde\chi^*(n-1,k),p_n} +
    \gen{p_n\tilde\chi^*(n-1,k-1), p_n}.
  \end{displaymath}
  If $k = 1$ we know from theorem \ref{teo:recursiveformula}
  that $H^1(M(d\scr A_{n-1}); \mat C) \cong H^1(M(d\scr A_{n-2}); \mat C)
  \oplus V_{(n-1,1)}$ and there is (at least) one copy of the standard representation at the
  degree $1$. Consider the case $k > 1$.
  By inductive hypothesis $\res{S_{n-1}}{S_n}\tilde\chi^*(n-1,k-1) = 
  \chi^*(n-1,k-1)$ contains exactly one copy of the standard representation
  therefore $\tilde\chi^*(n-1,k-1)$ must contain an irreducible representation
  which restricts to the standard representation of $S_{n-1}$;
  $V_{(n-1,1)}$ is not suitable because there is no copy of the trivial representation 
  in $\chi^*(n-1,k-1)$, so $\tilde\chi^*(n-1,k-1)$ must contain exactly one of the following
  \begin{displaymath}
    V_{(n-2,2)},\quad 
    V_{(n-2,1,1)}.
  \end{displaymath}
  
  Using Pieri's rule we see that both $V_{(n-2,1,1)} \otimes V_{(n-1,1)}$ and 
  $V_{(n-2,2)} \otimes V_{(n-1,1)}$ contain exactly
  one copy of the standard representation. 
  
  In particular $H^k(M(d\scr A_{n-1}); \mat C)$ contains
  exactly one copy of the standard representation 
  for every $1 < k < n-1$ and since there are $n-2$
  copies of the standard representation in $H^*(M(d\scr A_{n-1}); \mat C)$
  also $H^1(M(d\scr A_{n-1}); \mat C)$ contains exactly one copy
  of the standard representation.
\end{proof}

\begin{oss}\label{oss:quozienti}
  Proposition \ref{prop:wherearethestandard} can be used 
  for instance to compute the
  cohomology 
  of the quotient space $M(\scr A_{n-1})/S_{n-1}$. 
  Indeed, using the theorem on transfer
  we know that there is an isomorphism of graded algebras
  $H^*(M(\scr A_{n-1})/S_{n-1}; \mat C) \cong H^*(M(\scr A_{n-1}); \mat C)^{S_{n-1}}$.
  So, in order to compute the $\mat C$-vector space structure of
  $H^*(M(\scr A_{n-1})/S_{n-1}; \mat C)$ we need to look at those
  representations of $S_n$ whose restriction to $S_{n-1}$ contain
  a copy of the trivial representation, i.e. the trivial representation
  and the standard representation. 
  Therefore, when \(k=0\) or \(k=n-1\)   $H^k(M(\scr A_{n-1})/S_{n-1}; \mat C)$ is one dimensional,  
  while $H^k(M(\scr A_{n-1})/S_{n-1}; \mat C) $ is two dimensional when \(0<k<n-1\).
\end{oss}

\begin{prop}
  For $n \geq 3$ the copy of the standard tensor sign representation
  $V_{(2,1,\dots,1)}$ appears in the top cohomology $H^{n-2}(M(d\scr A_{n-1}); \mat C)$.
\end{prop}
\begin{proof}
  By induction on $n$; as in proposition \ref{prop:wherearethestandard}
  for $n = 3,4,5$ it follows from an explicit computation. 
  Let $n > 5$, from theorem \ref{teo:recursiveformula} we have
  \begin{displaymath}
    H^{n-2}(M(d\scr A_{n-1}); \mat C) \cong
    V_{(n-1,1)} \otimes H^{n-3}(M(d\scr A_{n-2}); \mat C).
  \end{displaymath}
  Again there must be exactly one irreducible representation of $S_n$ in
  $H^{n-3}(M(d\scr A_{n-2}); \mat C)$ whose restriction
  to $S_{n-1}$ contains a copy of $V_{(2,1,\dots,1)}$.
  This can't be $V_{(2,1,\dots,1)}$ because there
  is no copy of the alternating representation of
  $S_{n-1}$ in $H^{n-3}(M(d\scr A_{n-2}); \mat C)$.
  Therefore $H^{n-3}(M(d\scr A_{n-2}); \mat C)$ must
  contain one of the following representations
  of $S_n$:
  \begin{displaymath}
    V_{(2,2,1\dots,1)},\quad 
    V_{(3,1,\dots,1)}.
  \end{displaymath}
  But  $V_{(2,2,1,\dots,1)} \otimes V_{(n-1,1)}$
  and $V_{(3,1,\dots,1)} \otimes V_{(n-1,1)}$ contain exactly one copy of 
  $V_{(2,1,\dots,1)}$; therefore  $H^{n-2}(M(d\scr A_{n-1}); \mat C)$ contains exactly
  one copy of $V_{(2,1,\dots,1)}$.
\end{proof}

\subsection{The case $d$ odd}
From theorem \ref{teo:regular} we know that $H^*(C_n(d); \mat C)$
is the regular representation, in particular it contains 
$\dim V_{(n-1,1)} = n-1$ copies of the standard representation,
$\dim V_{(2,1, \dots, 1)} = n-1$ copies of the standard tensor
sign representation, one copy of the trivial and one copy of
the sign representations.

With the same argument as in proposition \ref{prop:wherearethestandard}
we can prove the following:
\begin{prop}\label{prop:wherearethestandardodd}
  For $n \geq 3$ and \(d\) odd there is exactly one copy of the 
  standard representation in the degree $k(d-1)$
  for each $1 \leq k \leq n-1$.
\end{prop}

\begin{oss}
  As in remark \ref{oss:quozienti}, proposition \ref{prop:wherearethestandardodd}
  can be used to compute the cohomology algebra of the quotient space
  $H^*(C_n(d)/S_{n-1}; \mat C)$. 
  In particular   $H^{k(d-1)}(C_n(d)/S_{n-1}; \mat C)$ is one dimensional  for every \(0\leq k \leq n-1\).
\end{oss}

Next we look at the sign representation; this was located by Lehrer
in \cite{lehrer2000ecc} using a formula for the generalized Poincar\`e polynomial. 
Our proof is different: we show an explicit generator.
\begin{prop}
  \label{prop:whereisthesign}
  Let $n = 2k$ or $n = 2k+1$ and \(d\) odd,
  then the copy of the sign representation appears
  in the component $H^{k(d-1)}(C_n(d); \mat C)$.
\end{prop}
\begin{proof}
  Consider the case $n = 2k$ and the following antisymmetrizer
  \begin{displaymath}
     x = \sum_{\sigma \in S_n} (-1)^\sigma A_{\sigma(1),\sigma(2)}A_{\sigma(3),\sigma(4)}\cdots 
    A_{\sigma(n-1),\sigma(n)} \in H^{k(d-1)}(C_n(d); \mat C).
  \end{displaymath}
  Of course $S_n$ acts on $\mat C x$ as $\tau x = (-1)^\tau x$,
  the non trivial part of  the argument consists in proving that $x \neq 0$.
  Consider the action of $S_n$ on the set of \emph{$2$-partitions of $\{1, \dots, n\}$}
  (that is partitions in which every block has cardinality $2$); 
  let $\Lambda$ be a $2$-partition and consider the following ordering on $\Lambda$
  \begin{displaymath}
    \Lambda = \{\Lambda_1, \dots, \Lambda_k\},\quad \Lambda_h = \{i_h, j_h\}
    \mbox{ with } i_h < j_h \mbox{ and } j_1 < \dots < j_k.
  \end{displaymath}
  In particular we can associate to every $\Lambda$ a permutation $\sigma_\Lambda \in S_n$
  such that $\sigma_\Lambda\{\{1,2\},\{3,4\}, \dots, \{n-1, n\}\} = \Lambda$ as follows
  \begin{displaymath}
    \sigma_\Lambda(2s) = j_s, \quad\quad  \sigma_\Lambda(2s+1) = i_{s+1}.
  \end{displaymath}
  Note that from this definition we have that $\sigma_\Lambda(A_{1,2}A_{3,4} \cdots A_{n-1,n})$
  is an element of the basis of admissible monomials (proposition \ref{prop:admissiblebasis}).

  Using the fact that $H^*(C_n(d); \mat C)$ is commutative and relation \(A_{i,j}=-A_{j,i}\)
  it can be easily seen that if $\tau \in S_n$ and 
  $\tau\{\{1,2\},\{3,4\}, \dots, \{n-1, n\}\} = \Lambda$  then
  \begin{displaymath}
    (-1)^\tau A_{\tau(1),\tau(2)} \cdots A_{\tau(n-1), \tau(n)} =
    (-1)^{\sigma_\Lambda} \sigma_\Lambda(A_{1,2}A_{3,4} \cdots A_{n-1,n}).
  \end{displaymath}
  In particular the expression of $x$ with respect to the basis of admissible
  monomials appears as follows
  \begin{equation}\label{eq:2partitions}
    x = m\sum_\Lambda (-1)^{\sigma_\Lambda} \sigma_\Lambda(A_{1,2}A_{3,4} \cdots A_{n-1,n})
  \end{equation}
  where $\Lambda$ runs over the $2$-partitions of $\{1, \dots, n\}$ and
  $m = k!2^k$ is the number of permutations of $S_n$ that fix the partition
  $\{\{1,2\}, \dots, \{n-1,n\}\}$, from which we conclude $x \neq 0$.
  
  Now consider the case $n = 2k+1$ and the element
  \begin{displaymath}
    x = \sum_{\sigma \in S_n} (-1)^\sigma A_{\sigma 1,\sigma 2}A_{\sigma 3,\sigma 4}\cdots 
    A_{\sigma n-2,\sigma n-1} \in H^{k(d-1)}(C_n(d); \mat C).
  \end{displaymath}
  With a similar argument as before we see that an analogous of \eqref{eq:2partitions}
  applies and therefore $x \neq 0$.
\end{proof}

Next we look at the standard tensor sign representation
$V_{(2,1, \dots, 1)}$.
\begin{prop}
  \label{prop:standardtensorsignconf}
  Consider \(d\) odd and \(k\geq 2\); if $n = 2k$  there is one copy of $V_{(2,1,\dots, 1)}$ in $H^{(k-1)(d-1)}(C_n(d); \mat C)$,
  one copy in $H^{k(d-1)}(C_n(d); \mat C)$,
  one copy in $H^{(n-1)(d-1)}(C_n(d); \mat C)$ and $2$ copies in each
  $H^{j(d-1)}(C_n(d); \mat C)$ for each $k < j < n-1$.
  If $n = 2k+1$ there is one copy of $V_{(2,1,\dots, 1)}$ in 
  $H^{k(d-1)}(C_n(d); \mat C)$, one copy in $H^{(n-1)(d-1)}(C_n(d); \mat C)$ 
  and $2$ copies in each $H^{j(d-1)}(C_n(d); \mat C)$ for each $k < j < n-1$.
\end{prop}
\begin{proof}
  By induction on $k$; the case  $k =2$  is trivial (see table \ref{table:decompositionodd}). 
  When $k > 2$, we use the recursive formula of
  theorem \ref{teo:recursiveformula}:
  \begin{displaymath}
    H^{j(d-1)}(C_n(d); \mat C) \cong H^{j(d-1)}(C_{n-1}(d); \mat C) \oplus
    \sx H^{(j-1)(d-1)}(C_{n-1}(d); \mat C) \otimes V_{(n-1,n)}\dx.
  \end{displaymath}
  Consider the case $n = 2k$.
  \begin{enumerate}[(i)]
  \item If $j = (k-1)$ then, by proposition
    \ref{prop:whereisthesign} we know that 
    the extended action on $H^{j(d-1)}(C_{n-1}(d); \mat C)$ must contain
    a copy of $V_{(2,1,\dots,1)}$.
  \item If $j = k$ then by inductive hypothesis
    the extended action of $S_n$ on $H^{(k-1)(d-1)}(C_{n-1}(d); \mat C)$ must
    contain an $S_n$-irreducible representation that restricts
    to $V_{(2,1,\dots,1)}$ and as in proposition
    \ref{prop:wherearethestandard} we know that
    $H^{k(d-1)}(C_n(d); \mat C)$ must contain a copy of
    $V_{(2,1,\dots,1)}$.
  \item If $j = (n-1)$ then $H^{(n-1)(d-1)}(C_n(d); \mat C)
    \cong H^{(n-2)(d-1)}(C_n(d); \mat C)\otimes V_{(n-1,1)}$ and
    as before $H^{(n-1)(d-1)}(C_n(d); \mat C)$ must contain a
    copy of $V_{(2,1,\dots,1)}$.
  \item if $k < j < n-1$ then $k-1 < j-1 < n-2$ and by inductive
    hypothesis the extended action on $H^{j(d-1)}(C_{n-1}(d); \mat C)$
    must contain two irreducible representations whose restrictions
    contain a copy of $V_{(2,1,\dots,1)}$; as before we conclude
    that $H^{j(d-1)}(C_n(d); \mat C)$ contains at least two copies
    of $V_{(2,1,\dots,1)}$.
  \end{enumerate}
  Observing that $H^*(C_n(d); \mat C) \cong \mat C S_n$ contains
  $n-1$ copies of $V_{(2,1,\dots,1)}$ we obtain the thesis.
  Now consider the case $n = 2k+1$.
  \begin{enumerate}[(i)]
  \item If $j = k$ then by inductive hypothesis
    we know that the $S_{n-1}$-action on 
    $H^{(k-1)(d-1)}(C_{n-1}(d); \mat C)$ contains a
    copy of $V_{(2,1,\dots,1)}$ and as before
    $H^{k(d-1)}(C_{n-1}(d); \mat C)$ contains a copy
    of $V_{(2,1,\dots,1)}$.
  \item If $j = k+1$, we know that the $S_{n-1}$ action
    on $H^{k(d-1)}(C_{n-1}(d); \mat C)$ contains
    a copy of the alternating representation and
    a copy of $V_{(2,1,\dots,1)}$. Anyway
    the extended action of $S_n$ on $H^{k(d-1)}(C_{n-1}(d); \mat C)$
    cannot contain a copy of $V_{(2,1,\dots, 1)}$ because
    $V_{(2,1,\dots, 1)} \otimes V_{(n-1,1)}$ contains a copy
    of the alternating representation (contradicting
    proposition \ref{prop:whereisthesign}). 
    Therefore
    the extended action on $H^{k(d-1)}(C_{n-1}(d); \mat C)$
    must contain a copy of the alternating representation of
    $S_n$ and an irreducible representation of $S_n$ whose
    restriction contains a copy of $V_{(2,1,\dots,1)}$.
    The copy of the alternating gives, after tensoring
    with $V_{(n-1,n)}$, one copy of $V_{(2,1,\dots,1)}$
    and the other irreducible representation gives
    another one.
  \item If $k+1 < j < n-1$ then $k < j-1 < n$ and
    as before we have that $H^{j(d-1)}(C_n(d); \mat C)$
    contains $2$ copies of $V_{(2,1,\dots,1)}$.
  \item If $j = n-1$ we have
    $H^{(n-1)(d-1)}(C_n(d); \mat C) \cong
    H^{(n-2)(d-1)}(C_{n-2}(d); \mat C) \otimes V_{(n-1,1)}$, which
    contains a copy of $V_{(2,1,\dots,1)}$.
  \end{enumerate}
  Again we conclude using the fact that $H^*(C_n(d); \mat C) \cong \mat C S_n$ contains
  $n-1$ copies of $V_{(2,1,\dots,1)}$.
\end{proof}

\section{The degrees $d-1$ and $2(d-1)$}
\label{sec:lowdegrees}
It is interesting to notice that  the recurrence formula of Theorem \ref{teo:recursiveformula} suffices to determine, 
for every \(n\geq 3\) and \(d\geq 2\),
an explicit decomposition of  $H^{d-1}(C_n(d);\mat C)$ and 
$H^{2(d-1)}(C_n(d); \mat C)$,    both as $S_n$ and as $S_{n+1}$-modules.

\subsection{The case $d$ even}
As in section \ref{sec:locatedeven} it suffices to study
the cohomology algebra of the \emph{deconed} braid arrangement
$H^*(M(d\scr A_{n-1}); \mat C)$; the isomorphism \eqref{eq:isodeconing} 
allows to infere formulas for the decomposition of $H^1(M(\scr A_{n-1}); \mat C)$ 
and $H^2(M(\scr A_{n-1}); \mat C)$ from the analogous formulas
for $M(d\scr A_{n-1})$.

\begin{prop}
  For every $n \geq 3$
  the following equality of $S_{n+1}$ modules
  holds:
  \begin{displaymath}
    H^1(M(d\scr A_{n-1}); \mat C) \cong V_{(n-1,2)}.
  \end{displaymath}
  In particular we have the following decomposition of $S_n$-modules:
  \begin{displaymath}
    H^1(M(d\scr A_{n-1}); \mat C) \cong V_{(n-1,1)} \oplus V_{(n-2,2)}.
  \end{displaymath}
\end{prop}
\begin{proof}
  By induction on $n$; we have already discussed the case $n = 3$ (see table \ref{table:decompositiondeconing}).
  Let $n > 3$, from theorem \ref{teo:recursiveformula} and 
  the inductive hypothesis we have
  \begin{displaymath}
    H^1(M(d\scr A_{n-1}); \mat C) \cong H^1(M(d\scr A_{n-1}); \mat C) \oplus
    V_{(n-1,1)} \cong V_{(n-2,2)} \oplus V_{(n-1,1)}
  \end{displaymath}
  and it is easily seen, using Pieri's rule, that $V_{(n-1,2)}$
  is the only representation of $S_{n+1}$ that restricts to
  $V_{(n-2,2)} \oplus V_{(n-1,1)}$.
\end{proof}

Next we look at $H^2(M(d\scr A_{n-1}); \mat C)$; its decomposition
can be recursively computed for $n \leq 6$ using theorem
\ref{teo:recursiveformula} and observing that for every 
$m < 6$ there exists a unique action of $S_{m+1}$ that
restricts to $H^2(M(d\scr A_{m-1}); \mat C)$ (see table \ref{table:decompositiondeconing}).
This way we obtain the following decomposition of $S_6$ modules:
\begin{align*}
  H^2(M(d\scr A_5); \mat C) \cong
  2\,{\tiny\yng(4,1,1)} \oplus 2\,{\tiny\yng(3,2,1)} \oplus
  {\tiny\yng(3,3)} \oplus {\tiny\yng(5,1)} \oplus
  {\tiny\yng(4,2)}.
\end{align*}
Again there is only one $S_7$-action that restricts
to $H^2(M(d\scr A_5); \mat C)$, namely
\begin{displaymath}
  H^2(M(d\scr A_5); \mat C) \cong {\tiny\yng(3,3,1)} \oplus {\tiny\yng(4,2,1)} \oplus
  {\tiny\yng(5,1,1)}.
\end{displaymath}

\begin{teo}
  For $n \geq 6$  the following equality of $S_{n+1}$-modules
  holds:
  \begin{displaymath}
    H^2(M(d\scr A_{n-1}); \mat C) \cong 
    V_{(n-1,1,1)} \oplus V_{(n-3,3,1)} \oplus V_{(n-2,2,1)}.
  \end{displaymath}
\end{teo}
\begin{proof}
  By induction on $n$; for $n = 6$ the result follows
  from our previous discussion. Let $n > 6$, from theorem \ref{teo:recursiveformula} and
  the inductive hypothesis we have
  \begin{displaymath}
    H^2(M(d\scr A_{n-1}); \mat C) \cong
    V_{(n-2,1,1)} \oplus V_{(n-4,3,1)} \oplus V_{(n-3,2,1)} \oplus
    \sx V_{(n-2,2)} \otimes V_{(n-1,1)}\dx
  \end{displaymath}
  Next we notice  that
  \begin{displaymath}
    V_{(n-2,2)} \otimes V_{(n-1,1)} \cong
    V_{(n-3,2,1)} \oplus V_{(n-3,3)} \oplus V_{(n-2,1,1)} \oplus
    V_{(n-2,2)} \oplus V_{(n-1,1)}.
  \end{displaymath}
  and therefore
  \begin{displaymath}
    H^2(M(d\scr A_{n-1}); \mat C) \cong
    2\,V_{(n-3,2,1)} \oplus V_{(n-3,3)} \oplus 2\,V_{(n-2,1,1)} \oplus
    V_{(n-2,2)} \oplus V_{(n-1,1)} \oplus V_{(n-4,3,1)}.
  \end{displaymath}

  Using Pieri's rule we see that the only irreducible representations of $S_{n+1}$
  whose restriction contains $V_{(n-3,2,1)}$ that can appear
  in the decomposition of the extended action on 
  $H^2(M(d\scr A_{n-1}); \mat C)$ are
  $V_{(n-3,3,1)}$ and $V_{(n-2,2,1)}$ and they must both 
  appear with multiplicity one. This forces the extended action
  of $S_{n+1}$ on $H^2(M(d\scr A_{n-1}); \mat C)$ to be
  \begin{displaymath}
    V_{(n-1,1,1)} \oplus V_{(n-3,3,1)} \oplus V_{(n-2,2,1)}.    
  \end{displaymath}
\end{proof}

\begin{oss}
  In particular for \(n\geq 7\)  we have the following decomposition of $S_n$-modules:
  \begin{align*}
    H^2(M(d\scr A_{n-1}); \mat C) \cong
    V_{(n-1,1)} \oplus 2\,V_{(n-2,1,1)} \oplus
    V_{(n-3,3)} \oplus 2\,V_{(n-3,2,1)} \oplus 
    V_{(n-4,3,1)} \oplus V_{(n-2,2)}
  \end{align*}
\end{oss}

\subsection{The case $d$ odd}
We have already discussed the decomposition of $H^{d-1}(C_n(d); \mat C)$
(proposition \ref{prop:decompositionhd-1odd}), so we only have to treat
the degree $2(d-1)$.

As before with an explicit computation it can be seen that 
\begin{displaymath}
  H^{2(d-1)}(C_5(d); \mat C) \cong
  {\tiny\yng(1,1,1,1,1)} \oplus {\tiny\yng(2,1,1,1)} \oplus
  {\tiny\yng(4,1)} \oplus 2\,{\tiny\yng(3,2)} \oplus 2\,\tiny\yng(2,2,1) \oplus
  {\tiny\yng(3,1,1)}.
\end{displaymath}
So, at first sight, there are two possible actions of $S_6$ that restrict
to $H^{2(d-1)}(C_5(d); \mat C)$, namely:
\begin{displaymath}
  {\tiny\yng(2,1,1,1,1)} \oplus 2\,{\tiny\yng(2,2,2)} \oplus
  {\tiny\yng(4,1,1)} \oplus 2\,{\tiny\yng(3,3)}\quad
  \mbox{ and }\quad
  {\tiny\yng(2,1,1,1,1)} \oplus {\tiny\yng(3,2,1)} \oplus
  {\tiny\yng(2,2,2)} \oplus {\tiny\yng(4,2)}.
\end{displaymath}
Anyway if the first case holds we would have
\begin{displaymath}
  H^{2(d-1)}(C_6(d); \mat C) \cong
  {\tiny\yng(2,1,1,1,1)} \oplus 2\,{\tiny\yng(2,2,2)} \oplus
  2\,{\tiny\yng(4,1,1)} \oplus 2\,{\tiny\yng(3,3)} \oplus
  {\tiny\yng(4,2)} \oplus {\tiny\yng(5,1)} \oplus
  {\tiny\yng(3,1,1,1)} \oplus {\tiny\yng(3,2,1)}
\end{displaymath}
which is not the restriction of an $S_7$ action;
therefore the second case must hold.

\begin{teo}
  For $n \geq 5$ and \(d\) odd 
  there is an isomorphism of $S_{n+1}$-modules
  \begin{displaymath}
    H^{2(d-1)}(C_n(d); \mat C) \cong
    V_{(n-3,1,1,1,1)} \oplus V_{(n-2,2,1)} \oplus
    V_{(n-3,2,2)} \oplus V_{(n-1,2)}.
  \end{displaymath}
\end{teo}
\begin{proof}
  First we observe that for every $n$ it holds:
  \begin{displaymath}
    V_{(n-2,1,1)} \otimes V_{(n-1,1)} \cong
    V_{(n-2,1,1)} \oplus V_{(n-3,2,1)} \oplus V_{(n-3,1,1,1)} \oplus
    V_{(n-1,1)} \oplus V_{(n-2,2)}.
  \end{displaymath}
  We prove the thesis by induction on $n$; we have already discussed
  the case $n = 5$. Let $n > 5$,
  from theorem \ref{teo:recursiveformula} and
  the inductive hypothesis we have
  \begin{align*}
    H^{2(d-1)}(C_n(d); \mat C) \cong H^{2(d-1)}(C_{n-1}(d); \mat C) \oplus
    \sx V_{(n-1,1)} \otimes H^{d-1}(C_{n-1}(d); \mat C)\dx \cong\\
    V_{(n-4,1,1,1,1)} \oplus 2\,V_{(n-3,2,1)} \oplus V_{(n-4,2,2)} \oplus
    2\,V_{(n-2,2)} \oplus V_{(n-2,1,1)} \oplus V_{(n-3,1,1,1)} \oplus V_{(n-1,1)}.
  \end{align*}
  The copy of $V_{(n-4,2,2)}$ can not appear as a component of the
  restriction of $V_{(n-4,2,2,1)}$ or $V_{(n-4,3,2)}$ (the latter makes sense
  only for $n \geq 7$) because there are no copies of
  $V_{(n-4,2,1,1)}$ and $V_{(n-4,3,1)}$ in $H^{2(d-1)}(C_n(d); \mat C)$.
  Therefore the extended action must contain a copy of
  $V_{(n-3,2,2)}$ and his restriction gives a copy
  of $V_{(n-4,2,2)}$ and a copy of $V_{(n-3,2,1)}$.
  The other copy of $V_{(n-3,2,1)}$ must appear as a component
  of the restriction of $V_{(n-2,2,1)}$ because there is only one copy
  of $V_{(n-4,2,2)}$ and there are no copies of $V_{(n-4,2,1,1)}$ and 
  $V_{(n-4,3,1)}$ in $H^{2(d-1)}(C_n(d); \mat C)$. The restriction of
  $V_{(n-2,2,1)}$ contains a copy of $V_{(n-2,2)}$, a copy of
  $V_{(n-2,1,1)}$ and a copy of $V_{(n-3,2,1)}$.
  Analogously the other copy of $V_{(n-2,2)}$ must appear
  as a component of the restriction of $V_{(n-1,2)}$;
  this gives a copy of $V_{(n-2,2)}$ and a copy of 
  $V_{(n-1,1)}$. At this point the copies of
  $V_{(n-4,1,1,1,1)}$ and $V_{(n-3,1,1,1)}$ must come
  from the restriction of $V_{(n-3,1,1,1,1)}$.

  Summarizing, there is only an action of $S_{n+1}$
  that restricts to the action $S_n \agisce
  H^{2(d-1)}(C_n(d); \mat C)$, namely:
  \begin{displaymath}
    V_{(n-3,1,1,1,1)} \oplus V_{(n-2,2,1)} \oplus
    V_{(n-3,2,2)} \oplus V_{(n-1,2)}.
  \end{displaymath}
\end{proof}

\begin{oss}
  In particular, for \(n\geq 6\) and \(d\) odd the following decomposition of $S_n$-modules
  holds:
  \begin{displaymath}
    H^{2(d-1)}(C_n(d); \mat C) \cong
    V_{(n-4,1,1,1,1)} \oplus V_{(n-3,1,1,1)} \oplus
    2\,V_{(n-3,2,1)} \oplus V_{(n-2,1,1)} \oplus 
    2\, V_{(n-2,2)} \oplus V_{(n-4,2,2)}  \oplus V_{(n-1,1)}.
  \end{displaymath}
\end{oss}

\bibliographystyle{plain}
\bibliography{arrangiamenti}

\begin{thebibliography}{10}

\bibitem{arnold1969crc}
V.I. Arnol'd.
\newblock {The cohomology ring of the colored braid group}.
\newblock {\em Mathematical Notes}, 5(2):138--140, 1969.

\bibitem{bjorner1992csc}
A.~Bjorner and G.M. Ziegler.
\newblock {Combinatorial stratification of complex arrangements}.
\newblock {\em Journal of the American Mathematical Society}, pages 105--149,
  1992.

\bibitem{cohen1976hil}
F.R. Cohen, T.J. Lada, and J.P. May.
\newblock {\em {The homology of iterated loop spaces}}.
\newblock Springer, 1976.

\bibitem{cohen1993representation}
F.R. Cohen and L.R. Taylor.
\newblock {On the representation theory associated to the cohomology of
  configuration spaces, from:``Algebraic topology (Oaxtepec, 1991)''}.
\newblock {\em Contemp. Math}, 146:91--109, 1993.

\bibitem{de1995wonderful}
C.~De~Concini and C.~Procesi.
\newblock {Wonderful models of subspace arrangements}.
\newblock {\em Selecta Mathematica, New Series}, 1(3):459--494, 1995.

\bibitem{delongueville2001crc}
M.~de~Longueville and C.A. Schultz.
\newblock {The cohomology rings of complements of subspace arrangements}.
\newblock {\em Mathematische Annalen}, 319(4):625--646, 2001.

\bibitem{deligne2000algebre}
P.~Deligne, M.~Goresky, and R.~MacPherson.
\newblock {L'algebre de cohomologie du complement, dans un espace affine, d'une
  famille finie de sous-espaces affines}.
\newblock {\em Michigan Math. J}, 48(1):121--136, 2000.

\bibitem{feichtner2000cac}
E.M. Feichtner and G.M. Ziegler.
\newblock {On cohomology algebras of complex subspace arrangements}.
\newblock {\em Transactions of the American Mathematical Society}, pages
  3523--3555, 2000.

\bibitem{gaiffi1996asn}
G.~Gaiffi.
\newblock {The actions of $S_{n+1}$ and $S_n$ on the cohomology ring of a
  Coxeter arrangement of type $A_{n-1}$}.
\newblock {\em manuscripta mathematica}, 91(1):83--94, 1996.

\bibitem{getzler1995ope}
E.~Getzler.
\newblock {Operads and moduli spaces of genus $0$ Riemann surfaces}.
\newblock In {\em {The moduli space of curves}}, volume 129 of {\em Prog.Math.}
  Birkh\"auser, Boston, 1995.

\bibitem{lehrer1987psa}
G.I. Lehrer.
\newblock {On the Poincar\'e series associated with Coxeter group actions on
  complements of hyperplanes}.
\newblock {\em J. London Math. Soc}, 36(2):275--294, 1987.

\bibitem{lehrer2000ecc}
G.I. Lehrer.
\newblock {Equivariant Cohomology of Configurations in $\mathbb{R}^d$}.
\newblock {\em Algebras and Representation Theory}, 3(4):377--384, 2000.

\bibitem{lehrer1986asg}
G.I. Lehrer and L.~Solomon.
\newblock {On the action of the symmetric group on the cohomology of the
  complement of its reflecting hyperplanes}.
\newblock {\em J. Algebra}, 104(2):410--424, 1986.

\bibitem{mathieu1996hsigman}
O.~Mathieu.
\newblock {Hidden $\Sigma_{n+1}$-Actions}.
\newblock {\em Commun. Math. Phys}, 176:467--474, 1996.

\bibitem{orlik1980cat}
P.~Orlik and L.~Solomon.
\newblock {Combinatorics and topology of complements of hyperplanes}.
\newblock {\em Inventiones Mathematicae}, 56(1):167--189, 1980.

\bibitem{orlik1992ah}
P.~Orlik and H.~Terao.
\newblock {\em {Arrangements of hyperplanes}}.
\newblock Springer, 1992.

\bibitem{robinson1996trs}
A.~Robinson and S.~Whitehouse.
\newblock {The tree representation of $S_{n+1}$}.
\newblock {\em Journal of Pure and Applied Algebra}, 111(1):245--254, 1996.

\bibitem{yuzvinsky2001osa}
S.~Yuzvinsky.
\newblock {Orlik-Solomon algebras in algebra and topology}.
\newblock {\em Russian Mathematical Surveys}, 56(2):293--364, 2001.

\end{thebibliography}

\end{document}